\pdfoutput=1
\RequirePackage{ifpdf}
\ifpdf 
\documentclass[pdftex]{sigma}
\else
\documentclass{sigma}
\fi

\numberwithin{equation}{section}

\newtheorem{Theorem}{Theorem}[section]
\newtheorem{Corollary}[Theorem]{Corollary}
\newtheorem{Lemma}[Theorem]{Lemma}
\newtheorem{Proposition}[Theorem]{Proposition}
 { \theoremstyle{definition}
\newtheorem{Definition}[Theorem]{Definition}
\newtheorem{Example}[Theorem]{Example}
\newtheorem{Remark}[Theorem]{Remark} }

\newcommand{\tolabel}[1]{\stackrel{#1}{\to}}

\renewcommand{\top}{\mathrm{top}}
\newcommand{\cE}{\mathcal{E}}
\newcommand{\cC}{\mathcal{C}}
\newcommand{\id}{\mathrm{id}}

\newcommand{\Rx}{\mathbb{R}^\times}
\newcommand{\R}{\mathbb{R}}
\newcommand{\Z}{\mathbb{Z}}

\newcommand{\ddu}[1]{\frac{d}{du}|_{u=0} \left[ #1 \right]}

\DeclareMathOperator{\Ber}{Ber}
\DeclareMathOperator{\im}{im}

\begin{document}

\allowdisplaybreaks

\renewcommand{\thefootnote}{$\star$}

\newcommand{\arXivNumber}{1502.06253}

\renewcommand{\PaperNumber}{058}

\FirstPageHeading

\ShortArticleName{Modular Classes of Lie Groupoid Representations up to Homotopy}

\ArticleName{Modular Classes of Lie Groupoid Representations\\ up to Homotopy\footnote{This paper is a~contribution to the Special Issue
on Poisson Geometry in Mathematics and Physics.
The full collection is available at \href{http://www.emis.de/journals/SIGMA/Poisson2014.html}{http://www.emis.de/journals/SIGMA/Poisson2014.html}}}

\Author{Rajan Amit MEHTA}

\AuthorNameForHeading{R.A.~Mehta}

\Address{Department of Mathematics \& Statistics, Smith College,\\ 44 College Lane, Northampton, MA 01063, USA}
\Email{\href{mailto:rmehta@smith.edu}{rmehta@smith.edu}}
\URLaddress{\url{http://math.smith.edu/~rmehta/}}

\ArticleDates{Received February 24, 2015, in f\/inal form July 23, 2015; Published online July 25, 2015}

\Abstract{We describe a construction of the modular class associated to a representation up to homotopy of a Lie groupoid. In the case of the adjoint representation up to homotopy, this class is the obstruction to the existence of a volume form, in the sense of Weinstein's ``The volume of a dif\/ferentiable stack''.}


\Classification{22A22; 53D17}

\renewcommand{\thefootnote}{\arabic{footnote}}
\setcounter{footnote}{0}

\section{Introduction}

The \emph{modular class} of a Lie groupoid was introduced by Evens, Lu, and Weinstein~\cite[Appendix~B]{elw}, simultaneously generalizing the modular character of a Lie group, which represents the failure of the Haar measure to be bi-invariant, and the modular class of a foliation~\cite{yamagami}. In light of Weinstein's notion of volume of a dif\/ferentiable stack~\cite{weinstein:volume}, the modular class of a~Lie groupoid can be interpreted as the obstruction to the existence of a volume form.

The Evens--Lu--Weinstein approach uses the fact that, if $G \rightrightarrows M$ is a Lie groupoid with Lie algebroid~$A$, then there is a natural representation of the $1$-jet prolongation groupoid $J^1 G$ on the $2$-term complex $A \to TM$. Although this representation does not unambiguously descend to a representation of $G$, the ambiguities are chain homotopic to $0$, so there is a well-def\/ined representation of~$G$ ``up to homotopy''. As a result, one can then obtain a true representation of $G$ on the ``determinant line bundle'' associated to the complex.

In recent years, a general def\/inition of representation up to homotopy of a Lie groupoid has been given by Arias Abad and Crainic \cite{aba-cra:rephomgpd}, but there are some dif\/f\/iculties that emerge when trying to carry over the Evens--Lu--Weinstein construction to this setting. One issue is that, in general, there is no analogue of the $1$-jet prolongation groupoid for an arbitrary representation up to homotopy (however, such an analogue does exist in the $2$-term case; see \cite{gra-me:vbgpd}). Another issue is that, although Abad-Crainic representations up to homotopy include the information of a (pseudo-)action of the Lie groupoid on a complex, the action is allowed to be degenerate, potentially leading to an ill-behaved action on the determinant line bundle.

\looseness=-1
In this paper, we show that the latter issue can be resolved. Although pseudo-actions are allowed to be degenerate at the level of chain maps, the axioms of a representation up to homotopy imply that the pseudo-actions are homotopy equivalent to nondegenerate chain maps. This allows us to construct a well-def\/ined true representation on the \emph{Berezinian line bundle}, which is the correct analogue of the determinant line bundle in this setting. Using this representation, we then construct the modular class. The modular class associated to a representation up to homotopy of a Lie groupoid~$G$ on a complex $(\cE, \partial)$ can be geometrically interpreted as the obstruction to the existence of a $G$-invariant volume element (in an appropriate homological sense) on~$\cE$.

Modular classes associated to Lie algebroids were def\/ined in~\cite{elw}, and since then a number of papers have appeared describing various constructions of modular classes associated to Lie algebroids (see~\cite{ks:modular} and the references therein). The relationship between the modular class of a~Lie groupoid and that of its associated Lie algebroid is given by the van Est map~\cite{crainic:vanest, mythesis, wx}.

It is known \cite{crainic:vanest, cf, fernandes, gs-m:vbalg, me:modules} that the modular class is only the f\/irst of a sequence of ``secondary'' characteristic classes associated to Lie algebroid representations up to homotopy. It would be interesting to f\/ind analogous constructions of higher characteristic classes for Lie groupoids. It would also be interesting to see whether the modular class is related to the existence of a metric, in the sense of del Hoyo and Fernandes~\cite{fernandes-hoyo}.

The structure of the paper is as follows. In Section~\ref{sec:cohomology}, we give a brief review of Lie groupoid cohomology. In Section~\ref{sec:line}, we def\/ine modular classes associated to line bundle representations of Lie groupoids. In Section~\ref{sec:vanest}, we show that the van Est map takes Lie groupoid modular classes to Lie algebroid modular classes. In Section~\ref{sec:modrep}, we explain how modular classes are def\/ined for Lie groupoid representations. The main results of the paper are in Section~\ref{sec:modrephom}, where modular classes are constructed for Lie groupoid representations up to homotopy. In Section~\ref{sec:adjoint}, we consider the case of the adjoint representation up to homotopy of a Lie groupoid and interpret the modular class as the obstruction to the existence of a volume form in the sense of~\cite{weinstein:volume}.

\section{Groupoid cohomology}\label{sec:cohomology}

Let $G \rightrightarrows M$ be a Lie groupoid. We will denote the source and target maps as $s, t\colon G \to M$. For $x \in M$, we denote the corresponding unit element by $1_x \in G$.

\looseness=-1
It is well-known that the (smooth) Eilenberg--Maclane cohomology of a Lie group extends to Lie groupoids in the following way. Let $G^{(k)}$ denote the space of composable $k$-tuples of elements of~$G$:
\begin{gather*} G^{(k)} = \big\{ (g_1,\dots, g_k) \in G^k \,|\,  s(g_i) = t(g_{i+1}) \big\},\end{gather*}
with $G^{(0)} = M$. The $k$-cochains are def\/ined as smooth functions on $G^{(k)}$.

The coboundary operator $\delta$ is given by $\delta f (g) = f(s(g)) - f(t(g))$ for $f \in C^\infty(M)$, and
\begin{gather*}
  \delta f(g_0, \dots, g_k) =   f(g_1, \dots, g_k) + \sum_{i=1}^{k-1} (-1)^i f(g_1, \dots, g_i g_{i+1}, \dots, g_k)
+ (-1)^k f(g_1, \dots, g_{k-1})
\end{gather*}
for $f \in C^\infty(G^{(k)})$, $k > 0$.

More generally, we can consider cochains with values in any abelian Lie group. For the present purposes, we will primarily be concerned with cochains taking values in the group~$\Rx$ of multiplicative real numbers, in which case one should replace addition and subtraction in the above formulas for~$\delta$ with multiplication and division.

\section{The characteristic class of a line bundle representation}
\label{sec:line}

Let $G \rightrightarrows M$ be a Lie groupoid, let $L \to M$ be a trivializable (real) line bundle, and let $\Delta$ be a representation of $G$ on $L$. That is, for each $g \in G$, we have a linear isomorphism $\Delta_g\colon L_{s(g)} \to L_{t(g)}$, smoothly depending on $g$, such that $\Delta_{gh} = \Delta_g \Delta_h$ for any composable pair of elements $(g,h) \in G^{(2)}$.

Choose a nonvanishing section $\sigma \in \Gamma(L)$. This allows us to def\/ine a function $\phi_\sigma \in C^\infty(G;\Rx)$ by the formula
\begin{gather}\label{eqn:modfn}
   \Delta_g \sigma_{s(g)} =  \phi_\sigma (g) \sigma_{t(g)}.
\end{gather}
We refer to $\phi_\sigma$ as the \emph{characteristic function} associated to the triple $(L,\Delta,\sigma)$.

\begin{Proposition}\quad
\begin{enumerate}\itemsep=0pt
\item[$1.$] The characteristic function $\phi_\sigma$ is a $1$-cocycle.
\item[$2.$] The cohomology class of $\phi_\sigma$ is independent of the choice of $\sigma$, and every representative of this cohomology class arises from some choice of $\sigma$.
\end{enumerate}
\end{Proposition}
\begin{proof}
Let $(g,h) \in G^{(2)}$ be a composable pair. From \eqref{eqn:modfn}, we have
\begin{gather*}
\Delta_h \sigma_{s(h)} = \phi_\sigma(h) \sigma_{t(h)}.
\end{gather*}
Applying $\Delta_g$ to both sides, we see that
\begin{gather*}
\Delta_{gh} \sigma_{s(h)} = \phi_\sigma(h) \phi_\sigma(g) \sigma_{t(g)}.
\end{gather*}
On the other hand, by directly putting $gh$ in for $g$ in \eqref{eqn:modfn}, we have
\begin{gather*}
\Delta_{gh} \sigma_{s(h)} = \phi_\sigma(gh) \sigma_{t(g)}.
\end{gather*}
By comparing the last two equations, we conclude that
\begin{gather*}
\phi_\sigma(g) \phi_\sigma(h) = \phi_\sigma(gh),
\end{gather*}
which is the cocycle condition.

For the second part of the proposition, let $\sigma$ and $\sigma'$ be two nonvanishing sections of $L$. We may write $\sigma' = f\sigma$ for some $f \in C^\infty(M;\Rx)$. Then for any $g \in G$, we have
\begin{gather*}
\Delta_g \sigma'_{s(g)} = \phi_{\sigma'}(g) \sigma'_{t(g)},
\end{gather*}
so
\begin{gather*}
f(s(g)) \Delta_g \sigma_{s(g)} = \phi_{\sigma'}(g) f(t(g)) \sigma_{t(g)}.
\end{gather*}
Comparing this with \eqref{eqn:modfn}, we see that
\begin{gather*}
\phi_{\sigma'}(g) = \phi_\sigma(g) \frac{f(s(g))}{f(t(g))} = \phi_\sigma(g) \delta f.
\end{gather*}
Since $f$ can be arbitrary, it follows that every representative of the cohomology class can be obtained.
\end{proof}

\begin{Definition}
     The class $\Phi_L := [\phi_\sigma] \in H^1(G;\Rx)$ is called the \emph{characteristic class} associated to the pair $(L, \Delta)$.
\end{Definition}

The characteristic class $\Phi_L$ classif\/ies trivializable line bundle representations up to isomorphism.

Suppose that $G$ has two trivializable line bundle representations $(L_1, \Delta^1)$ and $(L_2, \Delta^2)$. Then the tensor product $L_1 \otimes L_2$ has an induced representation $\Delta$ of $G$, given by
\begin{gather*}
     \Delta_g (\ell_1 \otimes \ell_2) = \Delta^1_g \ell_1 \otimes \Delta^2_g \ell_2.
\end{gather*}
A straightforward calculation shows that $\Phi_{L_1 \otimes L_2} = \Phi_{L_1}\Phi_{L_2}$.

\begin{Remark}\label{rmk:nontrivializable}
To include nontrivializable line bundles, we may def\/ine a slightly weaker characteristic class $\Phi^+_L \in H^1(G;\R^+)$ by taking
\begin{gather*}
     \Phi^+_L := \sqrt{|\Phi_{L \otimes L}|}.
\end{gather*}
This is well-def\/ined, since $L \otimes L$ is trivializable.

If $G$ has connected source-f\/ibers, then the characteristic classes $\Phi_L$ are always positive, so there is no information loss in passing to $\Phi^+_L$. One can then see that, if $G$ is source-connected, then line bundle representations are classif\/ied by the pair $(w_1(L), \Phi^+_L)$, where $w_1(L)$ is the Steifel--Whitney class.
\end{Remark}

\begin{Remark}
In the above, we have only considered real line bundles. The construction would work equally well for trivializable complex line bundles, in which case $\Phi_L$ would take values in the multiplicative complex numbers. However, it is unclear how the construction would extend to a complex line bundle with nontrivial Chern class.
\end{Remark}

\section{Relation with the van Est map}\label{sec:vanest}

Let $G \rightrightarrows M$ be a Lie groupoid with Lie algebroid $A \to M$. Recall that the van Est map \cite{crainic:vanest, mythesis, wx} is a homomorphism from the complex of real-valued (normalized) smooth groupoid cochains to the complex $\wedge \Gamma(A^*)$ of Lie algebroid cochains. To apply the van Est map to $\Rx$-valued cochains, we ned to f\/irst take the logarithm.

We will only require the van Est map in the case of $1$-cochains, where there is the following simple formula. For any $\psi \in C^1(G;\Rx)$ satisfying the normalization condition $\phi(g) = 1$, the image $V\psi \in \Gamma(A^*)$ is given by
\begin{gather}\label{eqn:vanest}
     V\psi(X) = \bar{X}(\log \psi) = \bar{X}(\psi)
\end{gather}
for any $X \in \Gamma(A)$, where $\bar{X}$ is the associated vector f\/ield along the submanifold of units of $G$, tangent to the $s$-f\/ibers.

Now, let $\Delta$ be a representation of $G$ on a trivializable line bundle $L$, and let $\nabla\colon \Gamma(A) \times \Gamma(L) \to \Gamma(L)$ be the induced representation of $A$ on $L$. For any nonvanishing section $\sigma \in \Gamma(L)$, there is an induced Lie algebroid $1$-cocycle $\theta_\sigma \in \Gamma(A^*)$, given by (see~\cite{elw})
\begin{gather}\label{eqn:modularalgbd}
     \theta_\sigma(X) \sigma = \nabla_X \sigma
\end{gather}
for $X \in \Gamma(A)$. The cocycle $\theta_\sigma$ represents the characteristic class associated to the representa\-tion~$\nabla$, as def\/ined in~\cite{elw}.

\begin{Theorem}
     For any fixed choice of nonvanishing section $\sigma \in \Gamma(L)$, the van Est map~\eqref{eqn:vanest} sends $\phi_\sigma$ to $\theta_\sigma$.
\end{Theorem}

\begin{proof}
To simplify formulas, we use $\sigma$ to identify $L$ with $M \times \R$. Correspondingly, $\Gamma(L)$ is identif\/ied with $C^\infty(M)$ via the map $f \sigma \mapsto f$. Then~\eqref{eqn:modfn} may be rewritten as
\begin{gather}\label{eqn:deltatriv}
     \Delta_g(s(g), r) = (t(g), \phi_\sigma(g) \cdot r)
\end{gather}
for $g \in G$.

     In order to explicitly describe the induced representation $\nabla$ of $A$ on $L$, let us represent an arbitrary element $a \in A_x$ by a path $\gamma_u$ in $G$, def\/ined for $u \in [0,\epsilon)$, satisfying the following properties:
\begin{enumerate}\itemsep=0pt
     \item[1)] $\gamma_0 = 1_x$,
     \item[2)] $s(\gamma_u) = x$ for all $u$, and
     \item[3)] $\ddu{\gamma_u} = a$.
\end{enumerate}
In property~(3), we are viewing $a$ as a vector in $T_{1_x} G$ that is tangent to the $s$-f\/ibers. Then the induced representation $\nabla$ is given by
\begin{gather*}
\ddu{\Delta_{\gamma_u}(x,f(x))} = (\rho(a), - D_{\rho(a)}f + \nabla_a f)
\end{gather*}
for $f \in C^\infty(M) \cong \Gamma(L)$. Here, $D_{\rho(a)}f$ is the directional derivative of $f$ in the direction of the vector $\rho(a)$. Using \eqref{eqn:deltatriv}, we then have
\begin{gather*}
\ddu{\phi_\sigma(\gamma_u) \cdot f(x)} = - D_{\rho(a)}f + \nabla_a f,
\end{gather*}
so
\begin{gather}\label{eqn:nablapoint}
     \nabla_a f = D_{\rho(a)}f + (D_a \phi_\sigma) \cdot f(x).
\end{gather}
Applying \eqref{eqn:nablapoint} to an entire section $X \in \Gamma(A)$, we have
\begin{gather}\label{eqn:nablasection}
     \nabla_X f = \rho(X)(f) + \bar{X}(\phi_\sigma) \cdot f.
\end{gather}
Recalling that $\sigma$ is identif\/ied with the constant function $1$ under our identif\/ication of~$\Gamma(L)$ with~$C^\infty(M)$, we may now use \eqref{eqn:vanest}, \eqref{eqn:modularalgbd}, and \eqref{eqn:nablasection} to obtain the result
\begin{gather*}
     \theta_\sigma(X) = \nabla_X 1 = \bar{X}(\phi_\sigma) = V\phi_\sigma(X).\tag*{\qed}
\end{gather*}
\renewcommand{\qed}{}
\end{proof}

\section{The modular class of a Lie groupoid representation}\label{sec:modrep}

Let $G \rightrightarrows M$ \looseness=-1 be a Lie groupoid, let $E \to M$ be a vector bundle, and let~$\Delta$ be a representation of~$G$ on~$E$. For each $k \geq 1$, there is a naturally induced representation~$\wedge^k \Delta$ of~$G$ on~$\wedge^k E$, given by
\begin{gather*}
     \wedge^k \Delta_g (e_1 \wedge \cdots \wedge e_k) = \Delta_g e_1 \wedge \cdots \wedge \Delta_g e_k
\end{gather*}
for $e_1, \dots, e_k \in E_{s(g)}$. In particular, $\hat{\Delta} := \wedge^\top \Delta$ is a representation of $G$ on the determinant line bundle $\det(E):= \wedge^\top E$. For simplicity, we assume that $E$ is orientable, so that $\det(E)$ is trivializable; however, the construction can be extended to the nonorientable case following Remark \ref{rmk:nontrivializable}.

\begin{Definition}
     The \emph{modular class} $\Phi_E$ of the representation $\Delta$ is def\/ined to be the characteristic class $\Phi_{\det(E)}$ associated to the line bundle representation $\hat{\Delta}$. If $\Phi_E = 1$, then the representation is said to be \emph{unimodular}.
\end{Definition}

Let $\sigma$ be a determinant element on $E$, i.e., a nonvanishing section of~$\det(E)$. Then the characteristic function~$\phi_\sigma$ associated to the triple $(\det(E), \hat{\Delta}, \sigma)$ measures the failure of~$\Delta$ to preserve $\sigma$. Thus the representation is unimodular if and only if there exists a $G$-invariant determinant element on~$E$.

\begin{Example}
     In the case where $G$ is a Lie group, every $\Delta_g$ is an automorphism of a vector space, so the determinant is def\/ined independently of any choices. We then recover the notion of the modular character, which measures the failure of the representation to be volume-preserving. In particular, the adjoint representation allows us to associate a canonical modular character to $G$. The adjoint representation is unimodular if and only if there exists a bi-invariant volume form on $G$, as is the case, for example, when $G$ is compact or nilpotent.
\end{Example}

\section{The modular class of a representation up to homotopy}\label{sec:modrephom}

\subsection{Berezinian bundles}
Let $\cE = E^0 \oplus E^1$ be a f\/inite-dimensional $\Z_2$-graded vector bundle over $M$. The \emph{Berezinian bundle} of $\cE$ is def\/ined as
\begin{gather*}
     \Ber(\cE) := \wedge^\top E^0 \otimes \wedge^\top \big(E^1\big)^*.
\end{gather*}
This line \looseness=-1 bundle is the supergeometric generalization of the determinant bundle, in sense that a nonvanishing section of $\Ber(\cE)$ allows one to def\/ine the \emph{superdeterminant} or \emph{Berezinian} of an invertible map between dif\/ferent f\/ibers of $\cE$. Specif\/ically, given an invertible degree-preserving linear map $T\colon  \cE_x \to \cE_y$ for $x,y \in M$, we may extend $T$ to a map $\hat{T}\colon  \Ber(\cE)_x \to \Ber(\cE)_y$ given by
\begin{gather}\label{eqn:inducedber}
     \hat{T}(e_1 \wedge \cdots \wedge e_k \otimes \xi_1 \wedge \cdots \wedge \xi_\ell) = Te_1 \wedge \cdots \wedge T e_k \otimes (T^*)^{-1} \xi_1 \wedge \cdots \wedge (T^*)^{-1} \xi_\ell
\end{gather}
for $e_i \in E^0_x$ and $\xi_i \in E^{1*}_x$. If $\sigma$ is a nonvanishing section of $\Ber(\cE)$, then $\Ber_\sigma (T) \in \Rx$ is def\/ined by the equation
\begin{gather*}
     (\Ber_\sigma (T))\sigma_y = \hat{T}(\sigma_x).
\end{gather*}
If it is possible to write $\sigma = \sigma^0 \otimes (\sigma^1)^*$, where $\sigma^i$ is a nonvanishing section of $\wedge^\top E^i$ for $i \in \Z_2$, then we have
\begin{gather}\label{eqn:berdet}
     \Ber_\sigma (T) = \frac{\det_{\sigma^0} (T^0)}{\det_{\sigma^1} (T^1)},
\end{gather}
where $T^i\colon E^i_x \to E^i_y$ are the components of $T$.

\begin{Definition}
     A f\/inite-dimensional $\Z_2$-graded vector bundle $\cE \to M$ is called \emph{superorientable} if $\Ber(\cE)$ is trivializable.
\end{Definition}

\subsection{Representations up to weak homotopy}

Let $G \rightrightarrows M$ be a Lie groupoid, and let $\cE = \bigoplus E^i$ be a $\Z$-graded cochain complex of vector bundles over $M$ with coboundary operator $\partial^i\colon E^i \to E^{i+1}$.

\begin{Definition}\label{dfn:repweakhom}
     A \emph{representation up to weak homotopy} of $G$ on the complex $(\cE, \partial)$ is a smooth map that associates to each $g \in G$ a chain map $\Delta_g\colon (\cE_{s(g)}, \partial) \to (\cE_{t(g)}, \partial)$, such that $\Delta_g \Delta_h$ is chain homotopic to $\Delta_{gh}$ for all $(g,h) \in G^{(2)}$. A representation up to weak homotopy is \emph{unital} if $\Delta_{1_x} = \id$ for all $x \in M$.
\end{Definition}

\begin{Remark}
     The notion of representation up to weak homotopy is, as the name suggests, weaker than the notion of representation up to homotopy in~\cite{aba-cra:rephomgpd}. In Def\/inition \ref{dfn:repweakhom}, we simply require that a chain homotopy exist between~$\Delta_g \Delta_h$ and $\Delta_{gh}$, whereas the def\/inition in~\cite{aba-cra:rephomgpd} includes a map associating to each~$(g,h) \in G^{(2)}$
     a specif\/ic chain homotopy $\Omega_{g,h}$ from~$\Delta_{gh}$ to~$\Delta_g \Delta_h$, as well as higher homotopy maps, associating to each $(g_1, \dots, g_{p+1}) \in G^{(p+1)}$ a chain $p$-homotopy, satisfying a series of coherence relations.

There is a functor from the category of representations up to homotopy to that of representations up to weak homotopy, given by forgetting the homotopy maps. Thus, the construction below may be applied to representations up to homotopy, although it is obviously insensitive to information contained in the homotopy maps.
\end{Remark}

\subsection{The induced representation on the Berezinian bundle}\label{sec:modular}

Let $(\cE, \partial)$ be a bounded $\Z$-graded cochain complex of f\/inite-dimensional vector bundles over $M$, and let $\Delta$ be a unital representation up to weak homotopy of $G$ on $(\cE, \partial)$.
By $\Ber(\cE)$, we mean the Berezinian bundle of the $\Z_2$-graded vector bundle obtained by quotienting out the grading on $\cE$ modulo $2$.

We wish to prove that $\Delta$ naturally induces a representation of $G$ on $\Ber(\cE)$. The key obstacle that needs to be overcome is the fact that $\Delta_g$ might not be invertible for all $g$, in which case the formula in \eqref{eqn:inducedber} cannot be used. We will see that this dif\/f\/iculty can be circumvented, but we f\/irst require the following brief detour into homological superalgebra.

\begin{Lemma}\label{lemma:strictchain}
          Let $\cC_1$ and $\cC_2$ be bounded $\Z$-graded cochain complexes of finite-dimensional vector spaces that are isomorphic as graded vector spaces $($in other words, $\cC_1$ and $\cC_2$ have the same dimension in each degree$)$. Let $f\colon \cC_1 \to \cC_2$ be a~chain homotopy equivalence. Then there exists a chain isomorphism $\tilde{f}\colon \cC_1 \to \cC_2$ that is chain homotopy equivalent to~$f$.
\end{Lemma}

\begin{proof}
     Assume without loss of generality that $\cC_1$ and $\cC_2$ both vanish below degree $0$. For each~$i$, we may (noncanonically) make a decomposition $C_1^i = B_1^i \oplus H_1^i \oplus B_1 ^{i+1}$, where $B_1^i = \im \partial_1^{i-1}$ and $H_1^i = B_1^i/\ker \partial_1^i$. We make a similar decomposition for $C_2^i$. The fact that there exists a~chain homotopy equivalence between $\cC_1$ and $\cC_2$ implies that $\dim H_1^i = \dim H_2^i$ for all $i$. Since we are assuming that $\dim C_1^i = \dim C_2^i$ for all $i$, we then have that
\begin{gather*}
     \dim B_1^i + \dim B_1^{i+1} = \dim B_2^i + \dim B_2^{i+1}
\end{gather*}
for all $i$. Because $\cC_1$ and $\cC_2$ are assumed to vanish below degree $0$, we have that $\dim B_1^0 = \dim B_2^0 = 0$, and by induction we conclude that $\dim B_1^i = \dim B_2^i$ for all~$i$.

Next, we write $f^i$ in block form with respect to the above decompositions. Since $f$ is a chain map, it is block upper-triangular:
\begin{gather}\label{eqn:fblock}
     f^i = \begin{bmatrix} f^{B^i} & * & * \\ 0 & f^{H^i} & * \\ 0 & 0 & f^{B^{i+1}}
           \end{bmatrix}.
\end{gather}
Because of the dimension arguments above, all of the blocks along the diagonal are square. Since $f$ is a chain homotopy equivalence, the block $f^{H^i}$ is invertible. Next, we will show that the blocks $f^{B^i}$ can be made invertible by a chain homotopy, which will complete the proof.

For any set of linear maps $\phi^i\colon B^i_1 \to B^i_2$, we can construct maps $\Phi^i\colon C^i_1 \to C^{i-1}_2$, which, in terms of the above decompositions, are of the form
\begin{gather*}
 \Phi^i = \begin{bmatrix} 0&0&0\\0&0&0\\\phi^i&0&0 \end{bmatrix}.
\end{gather*}
The maps $\Phi^i$ form a chain homotopy from $f^i$ to
\begin{gather*}
 f^i + \partial \circ \Phi^i + \Phi^{i+1} \circ \partial = \begin{bmatrix}f^{B^i} + \phi^i & * & * \\ 0 & f^{H^i} & * \\ 0 & 0 & f^{B^{i+1}} + \phi^{i+1}
\end{bmatrix}.
\end{gather*}
The result follows from the fact that the maps $\phi^i$ can be chosen arbitrarily.
\end{proof}

\begin{Lemma}\label{lemma:berident}
Let $\cC$ be a bounded $\Z$-graded cochain complex of finite-dimensional vector spaces, and let $f$ be a chain automorphism of $\cC$. If $f$ is chain homotopic to the identity, then $\Ber(f) = \Ber(\id) = 1$.
\end{Lemma}

\begin{proof}
     As in the proof of Lemma \ref{lemma:strictchain}, choose a decomposition $C^i = B^i \oplus H^i \oplus B^{i+1}$ for each~$i$. If~$f$ is a chain automorphism, then it has block form~\eqref{eqn:fblock}, where the blocks along the diagonal are nondegenerate. The Berezinian (see~\eqref{eqn:berdet}) of~$f$ is
\begin{gather}\label{eqn:berhom}
     \Ber(f) = \prod_i \frac{\det\big(f^{2i}\big)}{\det\big(f^{2i+1}\big)} = \prod_i \frac{\det\big(f^{H^{2i}}\big)}{\det\big(f^{H^{2i+1}}\big)}.
\end{gather}
If $f$ is chain homotopic to the identity, then~$f^{H^i} = I$ for all~$i$, and it follows that $\Ber(f) = 1$.
\end{proof}

The following is an immediate consequence of Lemma \ref{lemma:berident}.
\begin{Corollary}\label{cor:berhomotopy}
     Let $\cC_1$ and $\cC_2$ be as in Lemma~{\rm \ref{lemma:strictchain}}, and let $f_1, f_2\colon \cC_1 \to \cC_2$ be homotopy equivalent chain isomorphisms. Then the induced maps $\hat{f}_1, \hat{f}_2\colon \Ber(\cC_1) \to \Ber(\cC_2)$ are equal to each other.
\end{Corollary}

We now return to the situation described at the beginning of Section~\ref{sec:modular}.
\begin{Theorem}
Let $(\cE, \partial)$ be a bounded $\Z$-graded cochain complex of finite-dimensional vector bundles over $M$, and let $\Delta$ be a unital representation up to weak homotopy of $G$ on $(\cE, \partial)$. Then there is a canonically induced representation $\hat{\Delta}$ of $G$ on $\Ber(\cE)$.
\end{Theorem}
\begin{proof}
For any $g \in G$, the unital property of~$\Delta$ implies that $\Delta_g$ and $\Delta_{g^{-1}}$ are homotopy inverses. Therefore $\Delta_g$ is a chain homotopy equivalence from~$\cE_{s(g)}$ to~$\cE_{t(g)}$. If $\Delta_g$ is not invertible, then by Lemma~\ref{lemma:strictchain}, it can be replaced by a homotopy equivalent map which is invertible. By Corollary~\ref{cor:berhomotopy}, the induced map~$\hat{\Delta}_g$ from $\Ber(\cE)_{s(g)}$ to $\Ber(\cE)_{t(g)}$ does not depend on the choice of invertible replacement.

Since $\Delta_g \Delta_h$ is chain homotopic to $\Delta_{gh}$ for all $(g,h) \in G^{(2)}$, it follows from Corollary \ref{cor:berhomotopy} that $\hat{\Delta}_g\hat{\Delta}_h=\hat{\Delta}_{gh}$.
\end{proof}

For simplicity, we assume in the following def\/inition that $\cE$ is superorientable, but the construction can be extended to the nonsuperorientable case following Remark \ref{rmk:nontrivializable}.
\begin{Definition}
     The \emph{modular class} $\Phi_\cE$ of the representation up to weak homotopy $\Delta$ is def\/ined to be the characteristic class $\Phi_{\Ber(\cE)}$ associated to the line bundle representation $\hat{\Delta}$.
\end{Definition}

\begin{Remark}
     Let $\sigma$ be a nonvanishing section of $\Ber(\cE)$. Then the associated characteristic function is the appropriately-def\/ined Berezinian function, whose value at $g \in G$ is $\Ber_\sigma(\Delta_g)$, as in \eqref{eqn:berdet}. By ``appropriately def\/ined'', we mean that one may need to choose an invertible replacement for~$\Delta_g$ before taking the Berezinian. The representation up to weak homotopy is unimodular if and only if there exists a~$G$-invariant Berezinian element on~$\cE$.
\end{Remark}

\begin{Example}[regular representations up to homotopy]
A representation up to (weak) homotopy is called \emph{regular} if the coboundary map $\partial^i$ has constant rank for each $i$. In this case, the cohomology of the complex $(\cE, \partial)$ consists of vector bundles $H^i$, and there is a naturally induced representation $\Delta^i$ on $H^i$ for each $i$.

From \eqref{eqn:berhom}, we have the following relationship between the modular class~$\Phi_\cE$ and the modular classes~$\Phi_{H^i}$:
\begin{gather*}
 \Phi_\cE = \prod_i \frac{\Phi_{H^{2i}}}{\Phi_{H^{2i+1}}}.
\end{gather*}
\end{Example}

\section{The adjoint representation and volume elements}\label{sec:adjoint}
Recall \cite{aba-cra:rephomgpd,gra-me:vbgpd} that the adjoint representation of a Lie groupoid $G \rightrightarrows M$ is a representation up to homotopy on the $2$-term complex $A \tolabel{\rho} TM$, where $A$ is the Lie algebroid of $G$ and $\rho$ is the anchor map. Here, we view $A$ as being in degree $0$ and $TM$ as being in degree $1$. The adjoint representation depends on the choice of a ``horizontal lift,'' but dif\/ferent choices lead to representations up to homotopy that are homotopy equivalent to each other. Thus, for any f\/ixed nonvanishing section $\sigma$ of $\Ber(TM \oplus A) = \wedge^\top A \otimes \wedge^\top T^*M$, the resulting modular function~$\phi_\sigma$ is independent of the choice of horizontal lift. We will refer to the class $[\phi_\sigma] \in H^1(G)$ as \emph{the modular class of~$G$}.

In~\cite{weinstein:volume}, A.~Weinstein showed that $G$-invariant nonvanishing sections of $\wedge^\top A \otimes \wedge^\top T^*M$ may be interpreted as volume forms on the dif\/ferentiable stack~$M // G$ presented by~$G$. Super\-orientability of the graded vector bundle~$A \oplus TM$ is clearly a necessary condition for the existence of such a volume form. Under the assumption of superorientability, the modular class of~$G$ may be then be interpreted as the remaining obstruction.

We remark that, in the case where superorientability does not hold, one could instead consider a weaker version of the modular class, as indicated by Remark~\ref{rmk:nontrivializable}. We would expect that this class should be interpreted as the obstruction to the existence of a density on~$M // G$.

\begin{Example}[proper Lie groupoids]
In \cite{crainic:vanest}, it was shown that, if~$G$ is a proper Lie groupoid, then its cohomology vanishes in degrees $\geq 1$. Since the modular class lives in~$H^1(G)$, it obviously vanishes when~$G$ is proper. Therefore, assuming superorientability, proper Lie groupoids always admit volume forms in the sense of~\cite{weinstein:volume}.
\end{Example}

\begin{Example}[fundamental groupoids]
Let $M$ be a manifold, and consider the fundamental groupoid $\Pi_1(M) \rightrightarrows M$. Its adjoint representation up to homotopy is on the complex $TM \tolabel{\id} TM$, which is superorientable and has trivial cohomology. As a result, the associated modular class vanishes, so fundamental groupoids always admit volume forms.
\end{Example}

\subsection*{Acknowledgements}

We thank Ping Xu for helpful comments on a draft of the paper. We also thank the referees for many useful suggestions.

\pdfbookmark[1]{References}{ref}
\LastPageEnding

\end{document}